\documentclass[11pt]{article}
\usepackage{amsfonts}
\usepackage{mathrsfs}
\usepackage{amsmath}
\usepackage{amssymb}
\usepackage[mathscr]{eucal}
\usepackage{graphicx}
\usepackage{epstopdf}
\renewcommand{\paragraph}{\roman{paragraph}}

\newtheorem{theorem}{\scshape \mdseries  Theorem}[section]
\newtheorem{lemma}[theorem]{\scshape \mdseries  Lemma}

\begin{document}

\title{\sf Maximum reciprocal degree resistance distance index of unicyclic graphs}
\author{\bf Gui-Dong Yu$^{1,2,}$\thanks {Corresponding author.
E-mail address: guidongy@163.com (G.D. Yu)}, Xing-Xing Li$^{1}$, Gai-Xiang Cai$^{1}$\\
 {\small \it $1.$ School of Mathematics and Computation Sciences,
Anqing Normal University, Anqing 246133, P. R. China}\\
{\small \it $2.$ Basic Department, Hefei preschool education college, Hefei
230013, P. R. China}}
\date{}
\maketitle \thispagestyle{empty}

\noindent {\bf Abstract:} The reciprocal degree resistance distance index of a connected graph $G$ is defined as $RDR(G)=\sum\limits_{\{u,v\}\subseteq V(G)}\frac {d_G(u)+d_G(v)}{r_G(u,v)}$,
 where $r_G(u,v)$ is the resistance distance between vertices $u$ and $v$ in $G$. Let $\mathscr {U}_n$ denote the set of unicyclic graphs with $n$ vertices.
We study the graph with maximum  reciprocal degree resistance
distance index among all graphs in $\mathscr {U}_n$, and characterize
the corresponding extremal graph.


\noindent {\bf MR Subject Classifications:} 05C12, 05C35, 92E10

\section{Introduction}

Let $G=(V,E)$ be a simple connected graph of order $n$
with vertex set $V=V(G)=\{v_1,v_2,\ldots,v_n\}$ and edge set
$E=E(G)$. For any $u\in V(G)$, $d_G(u)$ is the degree of vertex $u$, the distance between vertices of $u$ and $v$, denoted by $d_G(u,v)$,
is the length of a shortest path between them.
Topological indices are numbers associated with molecular structures which serve for
quantitative relationships between chemical structures and properties. The first such
index was published by Wiener\cite{30}, but the name topological index was invented by
Hosoya \cite{16}. Many of them are based on the graph distance \cite{38}, the vertex degree \cite{10}.
In addition, several graph invariants are based on both the vertex degree and the graph
distance \cite{7}.

One of the most
intensively studied topological indices is the Wiener index. The Wiener index was introduced by American
chemist H. Wiener in \cite{30}, defined as
$$W(G)=\sum\limits_{\{u,v\}\subseteq V(G)}{d_G(u,v)}.$$

Another distance-based graph invariant Harary index, has been introduced by Plav\v{s}i\'{c} et al.\cite{8} and independently by Ivanciuc et al.\cite{18}
in 1993 for the characterization of molecular graphs. The Harary index $H(G)$ of graph $G$ is defined as
   $$H(G)=\sum\limits_{\{u,v\}\subseteq V(G)}\frac 1{d_G(u,v)}.$$

   For more results related to Harary index, please refer to \cite{8,9,25,29,32,33,34,35,39,41}.

The resistance distance $r(u,v)$ (if more than one graphs are considered, we write $r_G(u,v)$
in order to avoid confusion) between vertices $u$ and $v$ in $G$ is defined as the effective resistance between the two nodes of
the electronic network obtained so that its nodes correspond the vertices of $G$ and each edges of $G$ is replaced by a resistor of
unit resistance, which is compared by the methods of the theory of resistive electrucal networks based on Ohm's and Kirchoff's laws.

   The Kirchhoff index $Kf(G)$ of a graph $G$ is defined as \cite{2,19}

   $$Kf(G)=\sum\limits_{\{u,v\}\subseteq V(G)}r_G(u,v).$$

   As a new structure-descriptor, Kirchoff index is well studied, see recent papers \cite{3,6,21,22,23,26,28,36,37,40,42,43}.
   In 2017 Chen et al \cite{4} introduced a new graph invariant reciprocal to  Kirchoff index, named Resistance-Harary index, as
   $$RH(G)=\sum\limits_{\{u,v\}\subseteq V(G)}\frac 1{r_G(u,v)}.$$

For more results related to Resistance-Harary index, please refer to \cite{31,4}.

The first and the second Zagreb indices are defined as
$$M_1=M_1(G)=\sum\limits_{u\in V(G)}{d_G(u)}^2=\sum\limits_{uv\in E(G)}(d_G(u)+d_G(v)),$$
and
$$M_2=M_2(G)=\sum\limits_{uv\in E(G)}d_G(u)d_G(v),$$
respectively. These are the oldest \cite{11,12} and best studied degree-based topological indices (see the reviews \cite{13}, recent papers
\cite{14,20}, and the references cited therein).

Dobrynin and Kochetova \cite{5} and Gutman \cite{15} independently proposed a
vertex-degree-weighted version of Wiener index called degree distance, which is
defined for a connected graph G as
$$DD(G)=\frac {1}2\sum\limits_{\{u,v\}\in V(G)}(d_G(u)+d_G(v))d_G(u,v).$$

The reciprocal degree distance\cite{1} is defined as
$$RDD(G) = \frac {1}2\sum\limits_{\{u,v\}\in V(G)}\frac {d_G(u)+d_G(v)}{d_G(u,v)}.$$
Hua and Zhang \cite{17} have obtained lower and upper bounds for the reciprocal
degree distance of graph in terms of other graph invariants. The chemical applications and mathematical properties of the reciprocal degree distance are
well studied in \cite{1,24}.

   Analogous to the relationship between degree distance and reciprocal degree distance, we introduce
here a new graph invariant based on both the
vertex degree and the graph distance, named the reciprocal degree resistance distance index, as

   $$RDR(G)=\sum\limits_{\{u,v\}\subseteq V(G)}\frac {d_G(u)+d_G(v)}{r_G(u,v)}.$$

   A unicyclic graph is a connected graph with $n$ vertices and $n$
edges.  Let $\mathscr {U}_n$
denote the set of unicyclic graphs with $n$ vertices.
   In this paper, we determine the graph with maximum reciprocal degree resistance distance index among all graph in $\mathscr {U}_n$,  and
characterize the corresponding extremal graph.

\section{Preliminaries}

In this section, we will introduce some useful lemmas and three
transformations.

Let $C_g=v_1v_2\cdots v_gv_1$ be the cycle on $g\geq3$ vertices. For any two vertices $v_i,v_j\in V(C_g)$ with $i<j$, by Ohm's law, one has
$$r_{C_g}(v_i,v_j)=\frac {(j-i)(g+i-j)}g$$
\begin{lemma}\cite{19}
Let $x$ be a cut vertex of a connected graph $G$ and let $a$ and $b$ be vertices occurring in different components which arise upon deletion
of $x$.Then
$$r_G(a,b)=r_G(a,x)+r_G(x,b)$$.
\end{lemma}

\textbf{2.1 Edge-lifting transformation}

Let $G_1$ and $G_2$ be two graphs with $n_1\geq2$ and $n_2\geq2$
vertices, respectively. If $G$ is the graph obtained from $G_1$ and
$G_2$ by adding an edge between a vertex $u_0$ of $G_1$ to a vertex
$v_0$ of $G_2$. And $G{'}$ is the graph obtained by identifying
$u_0$ of $G_1$ to a vertex $v_0$ of $G_2$ and adding a pendent edge
to $u_0(v_0)$. We say that $G{'}$ is obtained from $G$ by an
edge-lifting transformation at $e=\{u_0,v_0\},$ see Fig 1.
\begin{figure}[h]
\centering
\includegraphics{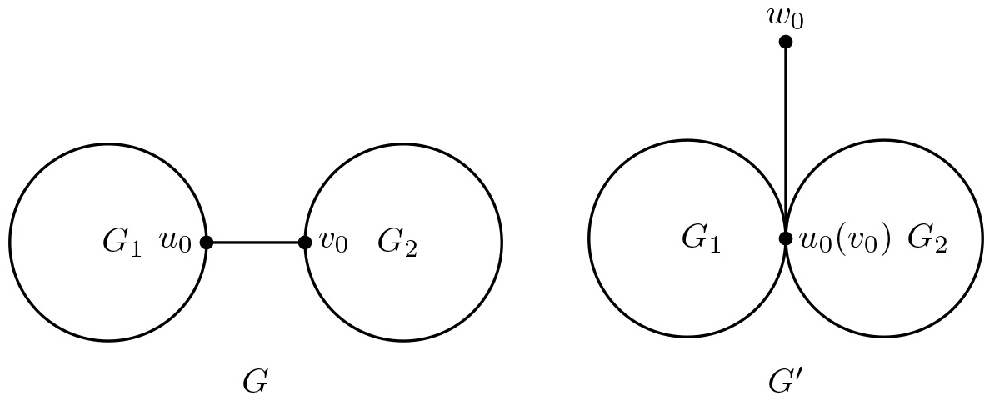}
\begin{center}
Fig 1: The edge-lifting transformation
\end{center}
\end{figure}

\begin{lemma}
If $G{'}$ can be obtained from $G$ by an edge-lifting transformation,
then $RDR(G)<RDR(G{'})$
\end{lemma}

{\bf Proof.} Let $G$ and $G{'}$ be shown in Fig.1.1. By the definition
of $RDR(G)$,
\begin{eqnarray*}
RDR(G)=&&\hspace{-.6cm}\sum\limits_{x,y\in V(G_1)\backslash\{u_0\}}\frac{d_G(x)+d_G(y)}{r_G(x,y)}+\sum\limits_{x,y\in V(G_2)\backslash\{v_0\}}\frac{d_G(x)+d_G(y)}{r_G(x,y)}\\
+&&\hspace{-.6cm}\sum\limits_{x\in V(G_1)\backslash\{u_0\}}\frac{d_G(u_0)+d_G(x)}{r_G(u_0,x)}+\sum\limits_{x\in V(G_2)\backslash\{v_0\}}\frac{d_G(u_0)+d_G(x)}{r_G(v_0,x)+1}\\
+&&\hspace{-.6cm}\sum\limits_{x\in V(G_1)\backslash\{u_0\}}\frac{d_G(v_0)+d_G(x)}{r_G(u_0,x)+1}
+\sum\limits_{x\in V(G_2)\backslash\{v_0\}}\frac{d_G(v_0)+d_G(x)}{r_G(v_0,x)}\\
+&&\hspace{-.6cm}\sum\limits_{x\in V(G_1)\backslash\{u_0\},y\in V(G_2)\backslash\{v_0\}}\frac{d_G(x)+d_G(y)}{r_G(x,u_0)+1+r_G(v_0,y)}+\frac{d_G(u_0)+d_G(v_0)}{r_G(u_0,v_0)}
\end{eqnarray*}

\begin{eqnarray*}
RDR(G{'})=&&\hspace{-.6cm}\sum\limits_{x,y\in V(G_1)\backslash\{u_0\}}\frac{d_{G{'}}(x)+d_{G{'}}(y)}{r_{G{'}}(x,y)}+\sum\limits_{x,y\in V(G_2)\backslash\{v_0\}}\frac{d_{G{'}}(x)+d_{G{'}}(y)}{r_{G{'}}(x,y)}\\
+&&\hspace{-.6cm}\sum\limits_{x\in V(G_1)\backslash\{u_0\}}\frac{d_{G{'}}(u_0)+d_{G{'}}(x)}{r_{G{'}}(u_0,x)}+\sum\limits_{x\in V(G_2)\backslash\{v_0\}}\frac{d_{G{'}}(v_0)+d_{G{'}}(x)}{r_{G{'}}(v_0,x)}\\
+&&\hspace{-.6cm}\sum\limits_{x\in
V(G_1)\backslash\{u_0\}}\frac{d_{G{'}}(w_0)+d_{G{'}}(x)}{r_{G{'}}(u_0,x)+1}
+\sum\limits_{x\in V(G_2)\backslash\{v_0\}}\frac{d_{G{'}}(w_0)+d_{G{'}}(x)}{r_{G{'}}(v_0,x)+1}\\
+&&\hspace{-.6cm}\sum\limits_{x\in V(G_1)\backslash\{u_0\},y\in
V(G_2)\backslash\{v_0\}}\frac{d_{G{'}}(x)+d_{G{'}}(y)}{r_{G{'}}(x,u_0)+r_{G{'}}(v_0,y)}+\frac{d_{G{'}}(u_0)+d_{G{'}}(w_0)}{r_{G{'}}(u_0,w_0)}
\end{eqnarray*}

(i) Noting that $d_G(x)+d_G(y)=d_{G{'}}(x)+d_{G{'}}(y)$,
$r_G(x,y)=r_{G{'}}(x,y)$ for $x,y\in V(G_1)\backslash\{u_0\}$ or
$x,y\in V(G_2)\backslash\{v_0\}$, and
$d_G(u_0)+d_G(v_0)=d_{G{'}}(u_0)+d_{G{'}}(w_0)$,
$r_G(u_0,v_0)=r_{G{'}}(u_0,w_0)$, then we have

$$\sum\limits_{x,y\in V(G_1)\backslash\{u_0\}}\frac{d_G(x)+d_G(y)}{r_G(x,y)}=\sum\limits_{x,y\in V(G_1)\backslash\{u_0\}}\frac{d_{G{'}}(x)+d_{G{'}}(y)}{r_{G{'}}(x,y)},$$

$$\sum\limits_{x,y\in V(G_2)\backslash\{v_0\}}\frac{d_G(x)+d_G(y)}{r_G(x,y)}=\sum\limits_{x,y\in V(G_2)\backslash\{v_0\}}\frac{d_{G{'}}(x)+d_{G{'}}(y)}{r_{G{'}}(x,y)},$$
$$\frac{d_G(u_0)+d_G(v_0)}{r_G(u_0,v_0)}=\frac{d_{G{'}}(u_0)+d_{G{'}}(w_0)}{r_{G{'}}(u_0,w_0)}.$$

(ii) Noting that $d_G(x)+d_G(y)=d_{G{'}}(x)+d_{G{'}}(y)$,
$r_G(x,u_0)+r_G(v_0,y)=r_{G{'}}(x,u_0)+r_{G{'}}(v_0,y)$ for $x\in
V(G_1)\backslash\{u_0\},~y\in V(G_2)\backslash\{v_0\}$, then we have
$$\sum\limits_{x\in V(G_1)\backslash\{u_0\},y\in
V(G_2)\backslash\{v_0\}}\frac{d_G(x)+d_G(y)}{r_G(x,u_0)+1+r_G(v_0,y)}<\sum\limits_{x\in
V(G_1)\backslash\{u_0\},y\in
V(G_2)\backslash\{u_0\}}\frac{d_{G{'}}(x)+d_{G{'}}(y)}{r_{G{'}}(x,u_0)+r_{G{'}}(v_0,y)}.$$

(iii) Noting that $r_{G{'}}(u_0,x)=r_G(u_0,x)$ for any $x\in
V(G_1)\backslash\{u_0\}$, and $d_{G{'}}(u_0)=d_G(u_0)+d_G(v_0)-1$,
then we have
$$(\sum\limits_{x\in V(G_1)\backslash\{u_0\}}\frac{d_{G{'}}(u_0)+d_{G{'}}(x)}{r_{G{'}}(u_0,x)}-\sum\limits_{x\in V(G_1)\backslash\{u_0\}}\frac{d_G(u_0)+d_G(x)}{r_G(u_0,x)})+$$
$$(\sum\limits_{x\in V(G_1)\backslash\{u_0\}}\frac{d_{G{'}}(w_0)+d_{G{'}}(x)}{r_{G{'}}(u_0,x)+1}-\sum\limits_{x\in V(G_1)\backslash\{u_0\}}\frac{d_G(v_0)+d_G(x)}{r_G(u_0,x)+1})$$
$$=\sum\limits_{x\in V(G_1)\backslash\{u_0\}}\frac{d_G(v_0)-1}{r_G(u_0,x)}+\sum\limits_{x\in V(G_1)\backslash\{u_0\}}\frac{1-d_G(v_0)}{r_G(u_0,x)+1}>0$$

(iv) Noting that $r_{G{'}}(v_0,x)=r_G(v_0,x)$, $d_{G{'}}(x)=d_G(x)$ for
any $x\in V(G_2)$, and $d_{G{'}}(v_0)=d_G(u_0)+d_G(v_0)-1$, then we
have
$$(\sum\limits_{x\in V(G_2)\backslash\{v_0\}}\frac{d_{G{'}}(v_0)+d_{G{'}}(x)}{r_{G{'}}(v_0,x)}-\sum\limits_{x\in V(G_2)\backslash\{v_0\}}\frac{d_G(v_0)+d_G(x)}{r_G(v_0,x)})+$$
$$(\sum\limits_{x\in V(G_2)\backslash\{v_0\}}\frac{d_{G{'}}(w_0)+d_{G{'}}(x)}{r_{G{'}}(v_0,x)+1}-\sum\limits_{x\in V(G_2)\backslash\{v_0\}}\frac{d_G(u_0)+d_G(x)}{r_G(v_0,x)+1})$$
$$=\sum\limits_{x\in V(G_2)\backslash\{v_0\}}\frac{d_G(u_0)-1}{r_G(v_0,x)}+\sum\limits_{x\in V(G_2)\backslash\{v_0\}}\frac{1-d_G(u_0)}{r_G(v_0,x)+1}>0.$$

Thus by (i)-(iv), we get $RDR({G{'}})-RDR(G)>0$, the results follows.
\hfill $\blacksquare$

\textbf{2.2 Cycle-lifting transformation}

Let $G$ be a graph as shown in Fig 2. Take a cycle $C$ in $G$,
say $C={v_1}{v_2}\cdots{v_p}{v_1},$  $G $ can be viewed as a graph
obtained by coalescing $C$ with a number of star subgraphs of $G$,
say $G_1,G_2,\cdots G_p$, by identifying $v_i$ with the center of
$G_i$ for all $i(1\leq i \leq p)$, denoted $|E(G_i)|=s_i$, and
$|V(G_i)|=s_i+1$. Deleting all edges in $G_i,$ joining $v_1$ to all
pendent vertices of $G_i(2\leq i \leq p)$, we obtained a new graph,
denoted by ${G{'}}$. (see Fig 2). This operation is called a
cycle-lifting transformation of $G$ with respect to $C$.
\begin{figure}[h]
\centering
\includegraphics{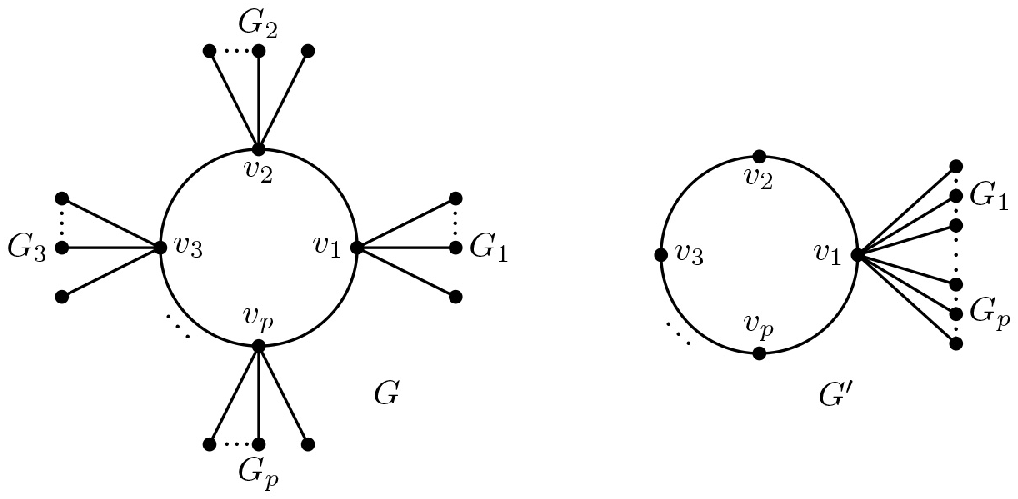}
\end{figure}
\begin{center}
Fig 2:   The cycle-lifting transformation
\end{center}
\begin{lemma}

If $G{'}$ can be obtained from $G$ by a cycle-lifting transformation, then $RDR(G)< RDR(G{'})$.
\end{lemma}

{\bf Proof.} Let $G$ and $G{'}$ be shown in Fig.1.2. By the definition of $RDR(G)$,
\begin{eqnarray*}
RDR(G)=&&\hspace{-.6cm}\sum_{i=1}^p\sum\limits_{x,y\in V(G_i)\backslash\{v_i\}}\frac{d_G(x)+d_G(y)}{r_G(x,y)}+\sum_{i=1}^p\sum\limits_{x\in V(G_i)\backslash\{v_i\}}\frac{d_G(v_i)+d_G(x)}{r_G(v_i,x)}\\
+&&\hspace{-.6cm}\sum_{\substack {i,j=1\\i\neq j}}^p\sum\limits_{\substack {x\in V(G_i)\backslash\{v_i\}\\y\in V(G_j)\backslash\{v_j\}}}\frac{d_G(x)+d_G(y)}{r_G(x,v_i)+r_G(v_i,v_j)+r_G(v_j,y)}+\sum_{i=1}^p\sum\limits_{\substack {x\in V(C_p)\backslash\{v_i\}\\y\in V(G_i)\backslash\{v_i\}}}\frac{d_G(x)+d_G(y)}{r_G(x,y)}\\
+&&\hspace{-.6cm}\sum\limits_{x,y\in V(C_p)}\frac{d_G(x)+d_G(y)}{r_G(x,y)}.
\end{eqnarray*}

\begin{eqnarray*}
RDR(G')=&&\hspace{-.6cm}\sum_{i=1}^p\sum\limits_{x,y\in V(G_i)\backslash\{v_1\}}\frac{d_{G{'}}(x)+d_{G{'}}(y)}{r_{G{'}}(x,y)}+\sum_{i=1}^p\sum\limits_{x\in V(G_i)\backslash\{v_1\}}\frac{d_{G{'}}(v_1)+d_{G{'}}(x)}{r_{G{'}}(v_1,x)}\\
+&&\hspace{-.6cm}\sum_{\substack {i,j=1\\i\neq j}}^p\sum\limits_{\substack {x\in V(G_i)\backslash\{v_1\}\\y\in V(G_j)\backslash\{v_1\}}}\frac{d_{G{'}}(x)+d_{G{'}}(y)}{r_{G{'}}(x,v_1)+r_{G{'}}(v_1,y)}+\sum_{i=1}^p\sum\limits_{\substack {x\in V(C_p)\backslash\{v_1\}\\y\in V(G_i)\backslash\{v_1\}}}\frac{d_{G{'}}(x)+d_{G{'}}(y)}{r_{G{'}}(x,y)}\\
+&&\hspace{-.6cm}\sum\limits_{x,y\in V(C_p)}\frac{d_{G{'}}(x)+d_{G{'}}(y)}{r_{G{'}}(x,y)}.
\end{eqnarray*}

(i) Noting that $d_G(x)+d_G(y)=d_{G{'}}(x)+d_{G{'}}(y)$,
$r_G(x,y)=r_{G{'}}(x,y)$ for $x,y\in V(G_i)\backslash\{v_i\}$ in $G$,  and $x,y\in V(G_i)\backslash\{v_1\}$ in ${G{'}}$, then we have

$$\sum_{i=1}^p\sum\limits_{x,y\in V(G_i)\backslash\{v_i\}}\frac{d_G(x)+d_G(y)}{r_G(x,y)}=\sum_{i=1}^p\sum\limits_{x,y\in V(G_i)\backslash\{v_1\}}\frac{d_{G{'}}(x)+d_{G{'}}(y)}{r_{G{'}}(x,y)}.$$

(ii) Let $r_G(v_i,v_j)=r_{{G{'}}}(v_i,v_j)=r_t$, when $|j-i|=t$, then $r_t=r_{p-t}$. $d_G(v_i)=s_i+2$ for any $v_i\in V(C_p)$ in $G$,  $d_{G{'}}(v_1)=s_1+\cdots+s_p+2$ for $v_1\in V(C_p)$ in $G{'}$, then we have
\begin{eqnarray*}
\sum\limits_{x,y\in V(C_p)}\frac{d_G(x)+d_G(y)}{r_G(x,y)}=\sum_{i=2}^p\frac{d_G(v_1)+d_G(v_i)}{r_G(v_1,v_i)}+\sum_{i=3}^p\frac{d_G(v_2)+d_G(v_i)}{r_G(v_2,v_i)}+\cdots +\frac{d_G(v_{p-1})+d_G(v_p)}{r_G(v_{p-1},v_p)}.
\end{eqnarray*}

Thus
\begin{eqnarray*}
\sum_{i=2}^p\frac{d_G(v_1)+d_G(v_i)}{r_G(v_1,v_i)}=&&\hspace{-.6cm}\frac{s_1+s_2+4}{r_1}+\frac{s_1+s_3+4}{r_2}+\cdots+\frac{s_1+s_p+4}{r_{p-1}},\\
\sum_{i=3}^p\frac{d_G(v_2)+d_G(v_i)}{r_G(v_2,v_i)}=&&\hspace{-.6cm}\frac{s_2+s_3+4}{r_1}+\cdots+\frac{s_2+s_p+4}{r_{p-2}},\\
&&\hspace{-.6cm}\cdots\\
\frac{d_G(v_{p-1})+d_G(v_p)}{r_G(v_{p-1},v_p)}=&&\hspace{-.6cm}\frac{s_{p-1}+s_p+4}{r_1}.\\
\end{eqnarray*}

Thus
\begin{eqnarray*}
\sum\limits_{x,y\in V(C_p)}\frac{d_G(x)+d_G(y)}{r_G(x,y)}=&&\hspace{-.6cm}\frac{(s_1+\cdots+s_{p-1})+(s_2+\cdots+s_p)+4(p-1)}{r_1}\\
+&&\hspace{-.6cm}\frac{(s_1+\cdots+s_{p-2})+(s_3+\cdots+s_p)+4(p-2)}{r_2}\\
+&&\hspace{-.6cm}\cdots+\frac{s_1+s_p+4}{r_{p-1}}.\\
\end{eqnarray*}

Similarly
\begin{eqnarray*}
\sum\limits_{x,y\in V(C_p)}\frac{d_{{G{'}}}(x)+d_{{G{'}}}(y)}{r_{{G{'}}}(x,y)}=&&\hspace{-.6cm}\frac{(s_1+\cdots+s_p)+4(p-1)}{r_1}+\frac{(s_1+\cdots+s_p)+4(p-2)}{r_2}\\
+&&\hspace{-.6cm}\cdots+\frac{s_1+\cdots+s_p+4}{r_{p-1}}.\\
\end{eqnarray*}

Then
\begin{eqnarray*}
\sum\limits_{x,y\in V(C_p)}\frac{d_{{G{'}}}(x)+d_{{G{'}}}(y)}{r_{G{'}}(x,y)}-\sum\limits_{x,y\in V(C_p)}\frac{d_G(x)+d_G(y)}{r_G(x,y)}=&&\hspace{-.6cm}\frac{s_2+\cdots+s_{p-1}}{r_1}+\frac{s_3+\cdots+s_{p-2}}{r_2}\\
+&&\hspace{-.6cm}\cdots+\frac{-(s_2+\cdots+s_{p-1})}{r_{p-1}}=0.
\end{eqnarray*}

(iii) Noting that $d_G(x)+d_G(y)=d_{G{'}}(x)+d_{G{'}}(y)$,
$r_G(x,v_i)+r_G(v_j,y)=r_{G{'}}(x,v_1)+r_{G{'}}(v_1,y)$ for $x\in
V(G_i)\backslash\{v_i\},~y\in V(G_j)\backslash\{v_j\}$ in $G$, and for $x\in
V(G_i)\backslash\{v_1\},~y\in V(G_j)\backslash\{v_1\}$ in $G{'}$, then we have
$$\sum_{\substack {i,j=1\\i\neq j}}^p\sum\limits_{\substack {x\in V(G_i)\backslash\{v_i\}\\y\in V(G_j)\backslash\{v_j\}}}\frac{d_G(x)+d_G(y)}{r_G(x,v_i)+r_G(v_i,v_j)+r_G(v_j,y)}<\sum_{\substack {i,j=1\\i\neq j}}^p\sum\limits_{\substack {x\in V(G_i)\backslash\{v_1\}\\y\in V(G_j)\backslash\{v_1\}}}\frac{d_{G{'}}(x)+d_{G{'}}(y)}{r_{G{'}}(x,v_1)+r_{G{'}}(v_1,y)}.$$

(iv) Noting that $r_G(v_i,x)=r_{G{'}}(v_1,x)=1$, $d_G(v_i)+d_G(x)=s_i+3$ for $x\in
V(G_i)\backslash\{v_i\}$ in $G$, $d_{G{'}}(v_1)+d_{G{'}}(x)=s_1+\cdots+s_p+3$ for $x\in
V(G_i)\backslash\{v_1\}$ in $G{'}$, then we have
$$\sum_{i=1}^p\sum\limits_{x\in V(G_i)\backslash\{v_i\}}\frac{d_G(v_i)+d_G(x)}{r_G(v_i,x)}=s_1(s_1+3)+s_2(s_2+3)+\cdots+s_p(s_p+3)
=\sum_{i=1}^ps_i(s_i+3).$$

$$\sum_{i=1}^p\sum\limits_{x\in V(G_i)\backslash\{v_1\}}\frac{d_{G{'}}(v_1)+d_{G{'}}(x)}{r_{G{'}}(v_1,x)}=s_1(s_1+\cdots+s_p+3)+\cdots+s_p(s_1+\cdots+s_p+3)=\sum_{i=1}^ps_i(s_1+\cdots+s_p+3).$$

Thus we have
$$\sum_{i=1}^p\sum\limits_{x\in V(G_i)\backslash\{v_1\}}\frac{d_{G{'}}(v_1)+d_{G{'}}(x)}{r_{G{'}}(v_1,x)}-\sum_{i=1}^p\sum\limits_{x\in V(G_i)\backslash\{v_i\}}\frac{d_G(v_i)+d_G(x)}{r_G(v_i,x)}=\sum_{\substack {i,j=1\\i\neq j}}^p {s_i\cdot s_j}.$$

(v) Noting that $d_G(x)+d_G(y)=3+s_i,d_{G{'}}(x)+d_{G{'}}(y)=3$,
$r_G(x,y)=r_{G{'}}(x,y)$ for any $x\in
V(C_p)\backslash\{v_i\},~y\in V(G_i)\backslash\{v_i\}$ in $G$, and for $x\in V(C_p)\backslash\{v_1\},~y\in V(G_i)\backslash\{v_1\}$ in $G{'}$, when $s_{i+t}>s_p(1\leq t\leq p-1),$ then $s_{i+t}=s_{i+t-p},$ then we have

\begin{eqnarray*}
\sum_{i=1}^p\sum\limits_{\substack {x\in V(C_p)\backslash\{v_i\}\\y\in V(G_i)\backslash\{v_i\}}}\frac{d_G(x)+d_G(y)}{r_G(x,y)}=&&\hspace{-.6cm}\sum\limits_{\substack {x\in V(C_p)\backslash\{v_1\}\\y\in V(G_1)\backslash\{v_1\}}}\frac{d_G(x)+d_G(y)}{r_G(x,y)}+\sum\limits_{\substack {x\in V(C_p)\backslash\{v_2\}\\y\in V(G_2)\backslash\{v_2\}}}\frac{d_G(x)+d_G(y)}{r_G(x,y)}\\
&&\hspace{-.6cm}+\cdots+\sum\limits_{\substack {x\in V(C_p)\backslash\{v_p\}\\y\in V(G_p)\backslash\{v_p\}}}\frac{d_G(x)+d_G(y)}{r_G(x,y)}.\\
\end{eqnarray*}

Thus
$$\sum\limits_{\substack {x\in V(C_p)\backslash\{v_1\}\\y\in V(G_1)\backslash\{v_1\}}}\frac{d_G(x)+d_G(y)}{r_G(x,y)}=\frac{s_1(3+s_2)}{r_1+1}+\frac{s_1(3+s_3)}{r_2+1}+\cdots+\frac{s_1(3+s_p)}{r_{p-1}+1},$$
$$\sum\limits_{\substack {x\in V(C_p)\backslash\{v_2\}\\y\in V(G_2)\backslash\{v_2\}}}\frac{d_G(x)+d_G(y)}{r_G(x,y)}=\frac{s_2(3+s_3)}{r_1+1}+\frac{s_2(3+s_4)}{r_2+1}+\cdots+\frac{s_2(3+s_1)}{r_{p-1}+1},$$
$$\cdots$$
$$\sum\limits_{\substack {x\in V(C_p)\backslash\{v_p\}\\y\in V(G_p)\backslash\{v_p\}}}\frac{d_G(x)+d_G(y)}{r_G(x,y)}=\frac{s_p(3+s_1)}{r_1+1}+\frac{s_p(3+s_2)}{r_2+1}+\cdots+\frac{s_p(3+s_{p-1})}{r_{p-1}+1},$$
thus we have
$$\sum_{i=1}^p\sum\limits_{\substack {x\in V(C_p)\backslash\{v_i\}\\y\in V(G_i)\backslash\{v_i\}}}\frac{d_G(x)+d_G(y)}{r_G(x,y)}=\sum_{i=1}^p\frac{s_i(3+s_{i+1})}{r_1+1}+\sum_{i=1}^p\frac{s_i(3+s_{i+2})}{r_2+1}+\cdots+\sum_{i=1}^p\frac{s_i(3+s_{i+p-1})}{r_{p-1}+1}.$$

Thus
$$\sum_{i=1}^p\sum\limits_{\substack {x\in V(C_p)\backslash\{v_i\}\\y\in V(G_i)\backslash\{v_i\}}}\frac{d_G(x)+d_G(y)}{r_G(x,y)}=\sum_{\substack {i,j=1\\i\neq j}}^p \frac{s_i\cdot(3+ s_j)}{r_G(v_i,v_j)+1}.$$

Similarly
\begin{eqnarray*}
\sum_{i=1}^p\sum\limits_{\substack {x\in V(C_p)\backslash\{v_1\}\\y\in V(G_i)\backslash\{v_1\}}}\frac{d_{G{'}}(x)+d_{G{'}}(y)}{r_{G{'}}(x,y)}=&&\hspace{-.6cm}\sum\limits_{\substack {x\in V(C_p)\backslash\{v_1\}\\y\in V(G_1)\backslash\{v_1\}}}\frac{d_{G{'}}(x)+d_{G{'}}(y)}{r_{G{'}}(x,y)}+\sum\limits_{\substack {x\in V(C_p)\backslash\{v_1\}\\y\in V(G_2)\backslash\{v_1\}}}\frac{d_{G{'}}(x)+d_{G{'}}(y)}{r_{G{'}}(x,y)}\\
&&\hspace{-.6cm}+\cdots+\sum\limits_{\substack {x\in V(C_p)\backslash\{v_1\}\\y\in V(G_p)\backslash\{v_1\}}}\frac{d_{G{'}}(x)+d_{G{'}}(y)}{r_{G{'}}(x,y)}.\\
\end{eqnarray*}
$$\sum\limits_{\substack {x\in V(C_p)\backslash\{v_1\}\\y\in V(G_1)\backslash\{v_1\}}}\frac{d_{G{'}}(x)+d_{G{'}}(y)}{r_{G{'}}(x,y)}=\frac{3s_1}{r_1+1}+\frac{3s_1}{r_2+1}+\cdots+\frac{3s_1}{r_{p-1}+1},$$
$$\sum\limits_{\substack {x\in V(C_p)\backslash\{v_1\}\\y\in V(G_2)\backslash\{v_1\}}}\frac{d_{G{'}}(x)+d_{G{'}}(y)}{r_{G{'}}(x,y)}=\frac{3s_2}{r_1+1}+\frac{3s_2}{r_2+1}+\cdots+\frac{3s_2}{r_{p-1}+1},$$
$$\cdots$$
$$\sum\limits_{\substack {x\in V(C_p)\backslash\{v_1\}\\y\in V(G_p)\backslash\{v_1\}}}\frac{d_{G{'}}(x)+d_{G{'}}(y)}{r_{G{'}}(x,y)}=\frac{3s_p}{r_1+1}+\frac{3s_p}{r_2+1}+\cdots+\frac{3s_p}{r_{p-1}+1}.$$

Thus we have
$$\sum_{i=1}^p\sum\limits_{\substack {x\in V(C_p)\backslash\{v_1\}\\y\in V(G_i)\backslash\{v_1\}}}\frac{d_{G{'}}(x)+d_{G{'}}(y)}{r_{G{'}}(x,y)}=\sum_{i=1}^p\frac{3s_i}{r_1+1}+\sum_{i=1}^p\frac{3s_i}{r_2+1}+\cdots+\sum_{i=1}^p\frac{3s_i}{r_{p-1}+1}.$$

Then
$$\sum_{i=1}^p\sum\limits_{\substack {x\in V(C_p)\backslash\{v_1\}\\y\in V(G_i)\backslash\{v_1\}}}\frac{d_{G{'}}(x)+d_{G{'}}(y)}{r_{G{'}}(x,y)}-\sum_{i=1}^p\sum\limits_{\substack {x\in V(C_p)\backslash\{v_i\}\\y\in V(G_i)\backslash\{v_i\}}}\frac{d_G(x)+d_G(y)}{r_G(x,y)}=-\sum_{\substack {i,j=1\\i\neq j}}^p \frac{s_i\cdot s_j}{r_G(v_i,v_j)+1}.$$

Thus by (i)-(v), we get $RDR({G{'}})-RDR(G)>0$, the results follows.
\hfill $\blacksquare$

\textbf{2.3 Cycle-shrinking transformation}

Denote by $S_n^p$ the unicyclic graph obtained from cycle $C_p$ by
attaching $n-p$ pendent edges to a vertex $v_1$ of $C_p$ (see
Fig 3). Let $G=S_n^p (p\geq4)$, deleting the edges
$v_{2}v_{3},\cdots, v_{p-1}v_{p}$, and adding the edges $v_{2}v_{p},
v_{3}v_{1}, \cdots,v_{p-1}v_{1}$, we obtain the graph $G'=S_n^3$.
This operation is called cycle-shrinking transformation. Denoted by
$W(W{'}$,resp), the set of pendent vertices of $G({G{'}}$,resp). Let
$|W|=k,$ $|W{'}|=k{'}$. It is clear that $p+k=n$, and $3+k{'}=n$.
Then, $k{'}=k+(p-3)$, thus $W{'}$ can be partitioned into two
subsets. One has $k$ vertices, which is naturally corresponding to
$W$. So we also denote it by $W$, another has $p-3$ vertices,
denoted by $\widetilde{W}$.
\begin{figure}[h]
\centering
\includegraphics{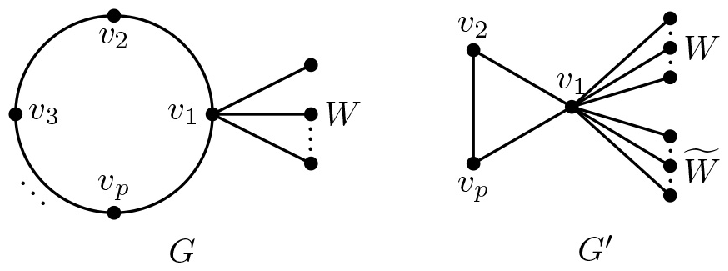}
\begin{center}
Fig 3:  The cycle-shrinking transformation
\end{center}
\end{figure}

\begin{lemma}
Let $G$ be an unicyclic graph of order $n$ and a cycle $C_p$ with $p\geq4$, if $G{'}$ can
be obtained from $G$ by a cycle-shrinking transformation, then
$RDR(G)< RDR(G{'}).$
\end{lemma}

 {\bf Proof.} Let $G$ and $G^{'}$ be shown in Fig.1.3. By the definition of $RDR(G)$,
\begin{eqnarray*}
RDR(G)=&&\hspace{-.6cm}\sum\limits_{x,y\in W}\frac{d_G(x)+d_G(y)}{r_G(x,y)}+\sum\limits_{x\in W}\frac{d_G(v_1)+d_G(x)}{r_G(v_1,x)}+\sum\limits_{\substack {x\in W\\y\in \{v_2,v_p\}}}\frac{d_G(x)+d_G(y)}{r_G(x,y)}\\
+&&\hspace{-.6cm}\sum\limits_{\substack {x\in W\\y\in V(C_p)\backslash\{v_1,v_2,v_p\}}}\frac{d_G(x)+d_G(y)}{r_G(x,y)}+\sum\limits_{x\in \{v_2,v_p\}}\frac{d_G(v_1)+d_G(x)}{r_G(v_1,x)}\\
+&&\hspace{-.6cm}\sum\limits_{x\in V(C_p)\backslash\{v_1,v_2,v_p\}}\frac{d_G(v_1)+d_G(x)}{r_G(v_1,x)}
+\sum\limits_{\substack {x\in V(C_p)\backslash\{v_1,v_2,v_p\}\\y\in \{v_2,v_p\}}}\frac{d_G(x)+d_G(y)}{r_G(x,y)}\\
+&&\hspace{-.6cm}\sum\limits_{x,y\in V(C_p)\backslash\{v_1,v_2,v_p\}}\frac{d_G(x)+d_G(y)}{r_G(x,y)}+\frac{d_G(v_2)+d_G(v_p)}{r_G(v_2,v_p)}.
\end{eqnarray*}

\begin{eqnarray*}
RDR(G')=&&\hspace{-.6cm}\sum\limits_{x,y\in W}\frac{d_{G{'}}(x)+d_{G{'}}(y)}{r_{G{'}}(x,y)}+\sum\limits_{x\in W}\frac{d_{G{'}}(v_1)+d_{G{'}}(x)}{r_{G{'}}(v_1,x)}+\sum\limits_{\substack {x\in W\\y\in \{v_2,v_p\}}}\frac{d_{G{'}}(x)+d_{G{'}}(y)}{r_{G{'}}(x,y)}\\
+&&\hspace{-.6cm}\sum\limits_{x\in W,y\in \widetilde{W}}\frac{d_{G{'}}(x)+d_{G{'}}(y)}{r_{G{'}}(x,y)}+\sum\limits_{x\in \{v_2,v_p\}}\frac{d_{G{'}}(v_1)+d_{G{'}}(x)}{r_{G{'}}(v_1,x)}+\sum\limits_{x\in \widetilde{W}}\frac{d_{G{'}}(v_1)+d_{G{'}}(x)}{r_{G{'}}(v_1,x)}\\
+&&\hspace{-.6cm}\sum\limits_{\substack {x\in \widetilde{W}\\y\in \{v_2,v_p\}}}\frac{d_{G{'}}(x)+d_{G{'}}(y)}{r_{G{'}}(x,y)}+\sum\limits_{x,y\in \widetilde{W}}\frac{d_{G{'}}(x)+d_{G{'}}(y)}{r_{G{'}}(x,y)}+\frac{d_{G{'}}(v_2)+d_{G{'}}(v_p)}{r_{G{'}}(v_2,v_p)}.
\end{eqnarray*}

(i) Noting that $d_G(x)+d_G(y)=d_{G{'}}(x)+d_{G{'}}(y)$, $r_G(x,y)=r_{G{'}}(x,y)$ for any $x,y\in W$, then we have
$$\sum\limits_{x,y\in W}\frac{d_G(x)+d_G(y)}{r_G(x,y)}=\sum\limits_{x,y\in W}\frac{d_{G{'}}(x)+d_{G{'}}(y)}{r_{G{'}}(x,y)}.$$

(ii) Noting that $d_G(x)+d_G(y)=d_{G{'}}(x)+d_{G{'}}(y)=3$, $r_G(x,y)=\frac {p-1}p+1\geq r_{G{'}}(x,y)=\frac {2}3+1$ for any $x\in W, y\in \{v_2,v_p\}$ in $G$ and $x\in W,y\in \{v_2,v_p\}$ in $G{'}$ , then we have
$$\sum\limits_{\substack {x\in W\\y\in \{v_2,v_p\}}}\frac{d_{G{'}}(x)+d_{G{'}}(y)}{r_{G{'}}(x,y)}-\sum\limits_{\substack {x\in W\\y\in \{v_2,v_p\}}}\frac{d_G(x)+d_G(y)}{r_G(x,y)}\geq 0.$$

(iii) Noting that $d_G(v_1)+d_G(x)=k+3$, $d_{G{'}}(v_1)+d_{G{'}}(x)=k+p$, $r_G(v_1,x)=r_{G{'}}(v_1,x)=1$, for $ x\in W$ in $G$ and $G{'}$ respectively, then we have
$$\sum\limits_{x\in W}\frac{d_{G{'}}(v_1)+d_{G{'}}(x)}{r_{G{'}}(v_1,x)}-\sum\limits_{x\in W}\frac{d_G(v_1)+d_G(x)}{r_G(v_1,x)}=k(k+p)-k(k+3)=k(p-3).$$

(iv) Noting that $d_G(x)+d_G(y)=3, r_G(x,y)=r_G(v_i,v_1)+1>2(3\leq i\leq p-1)$ for $x\in W,y\in V(C_p)\backslash\{v_1,v_2,v_p\}$, and $d_{G{'}}(x)+d_{G{'}}(y)=2,r_{G{'}}(x,y)=2$ for $x\in W,y\in \widetilde{W}$, then we have
$$\sum\limits_{\substack {x\in W\\y\in \widetilde{W}}}\frac{d_{G{'}}(x)+d_{G{'}}(y)}{r_{G{'}}(x,y)}-\sum\limits_{\substack {x\in W\\y\in V(C_p)\backslash\{v_1,v_2,v_p\}}}\frac{d_G(x)+d_G(y)}{r_G(x,y)}\geq k(p-3)-\frac{3k(p-3)}2=-\frac{k(p-3)}2.$$

Combing $(iii)-(iv)$, we have
$$(\sum\limits_{x\in W}\frac{d_{G{'}}(v_1)+d_{G{'}}(x)}{r_{G{'}}(v_1,x)}-\sum\limits_{x\in W}\frac{d_G(v_1)+d_G(x)}{r_G(v_1,x)})$$
$$+(\sum\limits_{\substack {x\in W\\y\in \widetilde{W}}}\frac{d_{G{'}}(x)+d_{G{'}}(y)}{r_{G{'}}(x,y)}-\sum\limits_{\substack {x\in W\\y\in V(C_p)\backslash\{v_1,v_2,v_p\}}}\frac{d_G(x)+d_G(y)}{r_G(x,y)})$$
$$\geq \frac{k(p-3)}2\geq 0.$$

(v) Noting that $d_G(v_1)+d_G(x)=k+4, r_G(v_1,x)=\frac {p-1}p$, and $d_{G{'}}(v_1)+d_{G{'}}(x)=k+p+1, r_{G{'}}(v_1,x)=\frac {2}3$, for $ x\in \{v_2,v_p\}$ in $G$ and $G{'}$ respectively, then we have
$$\sum\limits_{x\in \{v_2,v_p\}}\frac{d_{G{'}}(v_1)+d_{G{'}}(x)}{r_{G{'}}(v_1,x)}-\sum\limits_{x\in \{v_2,v_p\}}\frac{d_G(v_1)+d_G(x)}{r_G(v_1,x)}$$
$$=3(k+p+1)-\frac {2p(k+4)}{p-1}=\frac {(k+3p+1)(p-3)}{p-1}.$$

(vi) Noting that $d_G(v_1)+d_G(x)=k+4, r_G(v_1,x)=r_G(v_1,v_i)>1(3\leq i\leq p-1)$,for $x\in V(C_p)\backslash\{v_1,v_2,v_p\}$, and $d_{G{'}}(v_1)+d_{G{'}}(x)=k+p, r_{G{'}}(v_1,x)=1$, for $ x\in \widetilde{W}$, then we have
$$\sum\limits_{x\in \widetilde{W}}\frac{d_{G{'}}(v_1)+d_{G{'}}(x)}{r_{G{'}}(v_1,x)}-\sum\limits_{x\in V(C_p)\backslash\{v_1,v_2,v_p\}}\frac{d_G(v_1)+d_G(x)}{r_G(v_1,x)}$$
$$\geq (p-3)(k+p)-(p-3)(k+4)=(p-3)(p-4).$$

(vii) Noting that $d_G(x)+d_G(y)=4$£¬ for $x\in V(C_p)\backslash\{v_1,v_2,v_p\}$, $y\in \{v_2,v_p\}$, and $d_{G{'}}(x)+d_{G{'}}(y)=4$, $r_{G{'}}(x,y)=r_{G{'}}(v_2,v_1)+1=\frac 2{3}+1=\frac 5{3}$, for $x\in \widetilde{W}$, $y\in \{v_2,v_p\},$ by the symmetry between $v_2$ and $v_p$, then we have

\begin{eqnarray*}
&&\hspace{-.6cm}\sum\limits_{\substack {x\in \widetilde{W}\\y\in \{v_2,v_p\}}}\frac{d_{G{'}}(x)+d_{G{'}}(y)}{r_{G{'}}(x,y)}-\sum\limits_{\substack {x\in V(C_p)\backslash\{v_1,v_2,v_p\}\\y\in \{v_2,v_p\}}}\frac{d_G(x)+d_G(y)}{r_G(x,y)}\\
&&\hspace{-.6cm}=2\sum\limits_{i=3}^{p-1}\frac 3{r_{G{'}}(v_2,v_1)+1}-2\sum\limits_{i=3}^{p-1}\frac 4{r_G(v_2,v_i)}\\
&&\hspace{-.6cm}=\frac {18(p-3)}5-(\frac {8p}{p-1}+2\sum\limits_{i=4}^{p-1}\frac 4{r_G(v_2,v_i)}).\\
\end{eqnarray*}

Since $r_G(v_2,v_i)\geq r_G(v_2,v_4)=\frac {2(p-2)}p(4\leq i\leq p-1)$, if $\frac {2(p-2)}p\geq \frac 5{3}$, then $p\geq 12.$

Then
$$2\sum\limits_{i=4}^{p-1}\frac 4{r_G(v_2,v_i)})\leq \frac {24(p-4)}5$$

Thus

\begin{eqnarray*}
&&\hspace{-.6cm}\sum\limits_{\substack {x\in \widetilde{W}\\y\in \{v_2,v_p\}}}\frac{d_{G{'}}(x)+d_{G{'}}(y)}{r_{G{'}}(x,y)}-\sum\limits_{\substack {x\in V(C_p)\backslash\{v_1,v_2,v_p\}\\y\in \{v_2,v_p\}}}\frac{d_G(x)+d_G(y)}{r_G(x,y)}\\
&&\hspace{-.6cm}\geq \frac {18(p-3)}5-(\frac {8p}{p-1}+\frac {24(p-4)}5).\\
\end{eqnarray*}

(viii) Noting that $d_G(x)+d_G(y)=4$, $r_G(x,y)=r_G(v_i,v_j)(3\leq i<j\leq p-1)$, for $x,y\in V(C_p)\backslash\{v_1,v_2,v_p\}$, $d_{G{'}}(x)+d_{G{'}}(y)=2$, $r_{G{'}}(x,y)=2$, for $x,y\in \widetilde{W}$, then we have
\begin{eqnarray*}
&&\hspace{-.6cm}\sum\limits_{x,y\in \widetilde{W}}\frac{d_{G{'}}(x)+d_{G{'}}(y)}{r_{G{'}}(x,y)}-\sum\limits_{x,y\in V(C_p)\backslash\{v_1,v_2,v_p\}}\frac{d_G(x)+d_G(y)}{r_G(x,y)}\\
&&\hspace{-.6cm}=\frac {(p-3)(p-4)}2-\sum\limits_{\substack {i,j=3\\i<j}}^{p-1}\frac 4{r_G(v_i,v_j)}.\\
\end{eqnarray*}

Then
\begin{eqnarray*}
\sum\limits_{\substack {i,j=3\\i<j}}^{p-1}\frac 4{r_G(v_i,v_j)}=&&\hspace{-.6cm}\sum_{i=4}^p\frac{d_G(v_3)+d_G(v_i)}{r_G(v_3,v_i)}+\sum_{i=5}^p\frac{d_G(v_4)+d_G(v_i)}{r_G(v_4,v_i)}\\
&&\hspace{-.6cm}+\cdots +\frac{d_G(v_{p-2})+d_G(v_{p-1})}{r_G(v_{p-2},v_{p-1})}.\\
\end{eqnarray*}

When $j-i=t(i<j)$, $r_G(v_i,v_j)=r_t=r_{p-t}$,  thus
\begin{eqnarray*}
\sum_{i=4}^p\frac{d_G(v_3)+d_G(v_i)}{r_G(v_3,v_i)}=&&\hspace{-.6cm}\frac4{r_1}+\frac4{r_2}+\cdots+\frac4{r_{p-4}},\\
\sum_{i=5}^p\frac{d_G(v_4)+d_G(v_i)}{r_G(v_4,v_i)}=&&\hspace{-.6cm}\frac4{r_1}+\cdots+\frac4{r_{p-5}},\\
&&\hspace{-.6cm}\cdots\\
\frac{d_G(v_{p-2})+d_G(v_{p-1})}{r_G(v_{p-2},v_{p-1})}=&&\hspace{-.6cm}\frac4{r_1},\\
\end{eqnarray*}
thus
$$\sum\limits_{x,y\in V(C_p)\backslash\{v_1,v_2,v_p\}}\frac{d_G(x)+d_G(y)}{r_G(x,y)}=\frac{4(p-4)}{r_1}+\frac{4(p-5)}{r_2}+\cdots+\frac{4}{r_{p-4}}.$$

Since $r_1=\frac {p-1}p$, $r_t\geq r_2=\frac {2(p-2)}p(2\leq t\leq p-4)$, if $\frac {2(p-2)}p\geq \frac 5{3}$, then $p\geq 12.$

Then
$$\sum\limits_{x,y\in V(C_p)\backslash\{v_1,v_2,v_p\}}\frac{d_G(x)+d_G(y)}{r_G(x,y)}\leq \frac{4p(p-4)}{p-1}+\frac{6(p-4)(p-5)}5.$$

Thus
\begin{eqnarray*}
&&\hspace{-.6cm}\sum\limits_{x,y\in \widetilde{W}}\frac{d_{G{'}}(x)+d_{G{'}}(y)}{r_{G{'}}(x,y)}-\sum\limits_{x,y\in V(C_p)\backslash\{v_1,v_2,v_p\}}\frac{d_G(x)+d_G(y)}{r_G(x,y)}\\
&&\hspace{-.6cm}\geq\frac {(p-3)(p-4)}2-(\frac{4p(p-4)}{p-1}+\frac{6(p-4)(p-5)}5).\\
\end{eqnarray*}

When $p\geq12$, combing $(v)-(viii)$, we have
\begin{eqnarray*}
&&\hspace{-.6cm}(\sum\limits_{x\in \{v_2,v_p\}}\frac{d_{G{'}}(v_1)+d_{G{'}}(x)}{r_{G{'}}(v_1,x)}-\sum\limits_{x\in \{v_2,v_p\}}\frac{d_G(v_1)+d_G(x)}{r_G(v_1,x)})\\
+&&\hspace{-.6cm}(\sum\limits_{x\in \widetilde{W}}\frac{d_{G{'}}(v_1)+d_{G{'}}(x)}{r_{G{'}}(v_1,x)}-\sum\limits_{x\in V(C_p)\backslash\{v_1,v_2,v_p\}}\frac{d_G(v_1)+d_G(x)}{r_G(v_1,x)})\\
+&&\hspace{-.6cm}(\sum\limits_{\substack {x\in \widetilde{W}\\y\in \{v_2,v_p\}}}\frac{d_{G{'}}(x)+d_{G{'}}(y)}{r_{G{'}}(x,y)}-\sum\limits_{\substack {x\in V(C_p)\backslash\{v_1,v_2,v_p\}\\y\in \{v_2,v_p\}}}\frac{d_G(x)+d_G(y)}{r_G(x,y)})\\
+&&\hspace{-.6cm}(\sum\limits_{x,y\in \widetilde{W}}\frac{d_{G{'}}(x)+d_{G{'}}(y)}{r_{G{'}}(x,y)}-\sum\limits_{x,y\in V(C_p)\backslash\{v_1,v_2,v_p\}}\frac{d_G(x)+d_G(y)}{r_G(x,y)})\\
\geq &&\hspace{-.6cm}\frac {(k+3p+1)(p-3)}{p-1}+\frac {3(p-3)(p-4)}2+\frac {18(p-3)}5-\frac{4p(p-2)}{p-1}-\frac{6(p-1)(p-4)}5.\\
=&&\hspace{-.6cm}\frac {3p^3-22p^2+33p-54+k(10p-30)}{10(p-1)}.\\
\end{eqnarray*}

Let $f(p)=3p^3-22p^2+33p-54,  p\in [12,+\infty)$, this function is strictly increasing in this interval, thus $f(p)\geq f(12)=2358>0.$

Thus
$$\frac {3p^3-22p^2+33p-54+k(10p-30)}{10(p-1)}>0$$

(ix) Noting that $d_G(x)+d_G(y)=4$, $r_G
(v_2,v_p)=\frac {2(p-2)}p$ for $x\in \{v_2,v_p\}$ in $G$, and $d_{G{'}}(x)+d_{G{'}}(y)=4$, $r_{G{'}}(v_2,v_p)=\frac {2}3$
for $x\in \{v_2,v_p\}$ in $G{'}$, then we have
$$\frac{d_{G{'}}(v_2)+d_{G{'}}(v_p)}{r_{G{'}}(v_2,v_p)}-\frac{d_G(v_2)+d_G(v_p)}{r_G(v_2,v_p)}=6-\frac {2p}{p-2}=\frac {4p-12}{p-2}>0.$$

Thus by (i)-(ix), when $p\geq 12$, we get $RDR({G{'}})-RDR(G)>0$, the results
follows.

When $4\leq p\leq 11$, $U_n(G_0)=\{S_n^4,S_n^5,S_n^6,S_n^7,S_n^8,S_n^9,S_n^{10},S_n^{11}\}$, By comparison of direct calculation, we can obtained $RDR(S_n^3)>RDR(S_n^p)$.

The proof is completed.\hfill $\blacksquare$
\section {Maximum reciprocal degree resistance distance index of unicyclic graphs}

Let $G$ be a connected graph with exactly one cycle, say
$C_p=v_1v_2\cdots v_pv_1$. $G$ can be viewed as a graph obtained by
coalescing $C_p$ with a number of trees,
 denoted by $G_1,G_2,\cdots G_p$, by identifying $v_i$ with some vertex of $G_i$ for all $i(1\leq i\leq p)$. Denote $S_n^3$ by
the unicyclic graph obtained from cycle $C_3$ by attaching $n-3$ pendent edges to a vertex of $C_3$.  The next theorem discusses the maximum reciprocal degree resistance distance index in $\mathscr {U}_n$, which is the main result of our work.

\begin{theorem} Let $G$ be a unicyclic graph of order $n$, then $RDR(G)\leq RDR(S_n^3)$, with equality holds if and only if $G\cong S_n^3$.
\end{theorem}

{\bf Proof.} Let $G$ be a unicyclic graph of order $n$ such that the $\xi(G)$ is as big as possible.
By Lemma 2.2-2.4, we have $RDR(G)\leq RDR(S_n^3)$, equality holds if and only if $G\cong S_n^3$.\hfill $\blacksquare$

{\bf Conflicts of Interests} \noindent

 The authors declare that there is no conflict
of interests regarding the publication of this paper.

{\bf Acknowledgments} \noindent

This work is jointly supported by the National Natural Science
Foundation of China (11671164), the Natural Science Foundation of
Anhui Province (1808085MA04), and the Natural Science Foundation of
Department of Education of Anhui Province (KJ2017A362).

All authors are very grateful to Professor Shu-Chao Li of Center China Normal University. Because this paper is completed under the guidance of Prof. Li.

\end{document}